\documentclass[12pt,twoside]{amsart}
\usepackage{amssymb}
\usepackage{amscd}

\title[Kleiman-Mori cone]
{On the Kleiman-Mori cone}
\author{Osamu Fujino} 
\date{2005/1/5}
\address{Graduate School of Mathematics\\ 
 Nagoya University, Chikusa-ku Nagoya 464-8602 Japan}
\email{fujino@math.nagoya-u.ac.jp}
\newcommand{\Pic}[0]{{\operatorname{Pic}}}
\newcommand{\AD}[0]{{\operatorname{AD}}}

\newcommand{\Hom}[0]{{\operatorname{Hom}}}

\newcommand{\Proj}[0]{{\operatorname{Proj}}}
\newtheorem{thm}{Theorem}[section]
\newtheorem{lem}[thm]{Lemma}
\newtheorem{cor}[thm]{Corollary}
\newtheorem{prop}[thm]{Proposition}
\newtheorem{que}[thm]{Question}
\newtheorem{conj}[thm]{Conjecture}

\theoremstyle{definition}
\newtheorem{ex}[thm]{Example}
\newtheorem{defn}[thm]{Definition}
\newtheorem{rem}[thm]{Remark}
\newtheorem*{ack}{Acknowledgments}       
 
\newtheorem*{notation}{Notation}         
       
\newtheorem{say}[thm]{}
\begin{document}
\bibliographystyle{amsalpha+}

\begin{abstract}
The Kleiman-Mori cone plays important roles in the birational 
geometry. 
In this paper, we construct complete varieties whose 
Kleiman-Mori cones have interesting properties. 
First, we construct a simple and explicit example of complete 
non-projective singular varieties for which 
Kleiman's ampleness criterion does not hold. 
More precisely, we construct a complete non-projective toric variety $X$ 
and a line bundle $L$ on $X$ such that $L$ is positive 
on $\overline {NE}(X)\setminus \{0\}$. 
Next, we construct complete singular varieties $X$ with 
$NE(X)=N_1(X)\simeq \mathbb R^k$ for any $k$. 
These explicit examples seem to be missing in the literature. 
\end{abstract}

\maketitle

\section{{Introduction}}\label{intro}

The Kleiman-Mori cone plays important roles in the birational 
geometry. 
In this paper, we construct complete varieties whose 
Kleiman-Mori cones have interesting properties. 
First, we construct a simple and explicit example of 
complete non-projective 
singular varieties for which Kleiman's ampleness criterion 
does not hold. 
More precisely, we construct a complete non-projective toric variety $X$ 
and a line bundle $L$ on $X$ such that $L$ is positive 
on $\overline {NE}(X)\setminus \{0\}$. 
\begin{defn} 
Let $V$ be 
a complete algebraic scheme defined over an algebraically 
closed field $k$. 
We say that {\em{Kleiman's ampleness criterion holds}} for $V$ if and 
only if the interior of the nef cone of $V$ coincides with the ample 
cone of $V$. 
\end{defn} 
Note that Kleiman's original statements are very sharp. 
We recommend the readers to see \cite[Chapter IV \S2 Theorems 
1, 2]{kl}. 
Of course, our example is not \lq\lq quasi-divisorial\rq\rq 
in the sense of Kleiman (see \cite[Chapter IV \S2 
Definition 4 and Theorem 2]{kl}). 
We do not repeat the definition of {\em{quasi-divisorial}} since 
we do not use it in this paper. 
Note that if $X$ is projective or $\mathbb Q$-factorial then 
$X$ is quasi-divisorial in the sense of Kleiman. 
Next, we construct complete singular varieties $X$ with 
$NE(X)=N_1(X)\simeq \mathbb R^k$ for any $k$. We note 
that the condition $NE(X)=N_1(X)$ is equivalent to 
the following one:~a line bundle $L$ is nef 
if and only if $L$ is numerically trivial. 
These explicit examples seem to be missing in the literature. 
We adopt the toric geometry to construct examples. 

\begin{notation} 
We freely use the basic notation of the toric geometry 
throughout this paper. 
Let $v_i\in N$ for $0\leq i\leq k$. 
Then the symbol $\langle v_{1}, v_{2}, 
\cdots, v_{k}\rangle$ denotes the cone 
$\mathbb R_{\geq 0}v_{1}+\mathbb R_{\geq 0}v_2+\cdots 
+\mathbb R_{\geq 0}v_{k}$ in $N_{\mathbb R}$. 
\end{notation}

\begin{ack}
I would like to thank Professor 
Shigefumi Mori and 
Doctor Hiroshi Sato for fruitful discussions and 
useful comments. 
\end{ack}

\section{On Kleiman's ampleness criterion}

In this section, we construct explicit examples for 
which Kleiman's ampleness criterion does not hold. 
We think that the following example is the 
simplest one. It seems to be 
easy to construct a lot of singular toric varieties 
for which Kleiman's ampleness criterion does not 
hold. The reader can find many examples of 
singular toric $3$-folds in \cite{fs}. 
He can easily check that Kleiman's ampleness 
criterion does not hold for $X_6$ in \cite{fs}. 
For the cone theorem 
of toric varieties, see \cite[Theorem 4.1]{part1}. 
\begin{say}[Construction] 
We fix $N= \mathbb Z^3$. 
We put 
\begin{align*}
v_1  &= (1,0,1), & v_2&=(0,1,1), & v_3&=(-1,-1,1),\\
v_4  &= (1,0,-1), & v_5&=(0,1,-1), & v_6&=(-1,-1,-1).
\end{align*} 
We consider the following fans. 
$$
\Delta_P=
\left \{
\begin{array}{ccc}
\langle v_1, v_2, v_4\rangle, &
\langle v_2, v_4, v_5\rangle, &
\langle v_2, v_3, v_5, v_6\rangle, \\
\langle v_1, v_3, v_4, v_6\rangle, &
\langle v_1, v_2, v_3\rangle, &
\langle v_4, v_5, v_6\rangle, \\ 
\text{and their faces}& &
\end{array}
\right \}, $$ 
and 
$$
\Delta_Q=
\left\{ 
\begin{array}{ccc}
\langle v_1, v_2, v_4, v_5\rangle, &
\langle v_2, v_3, v_5, v_6\rangle, &
\langle v_1, v_3, v_4, v_6\rangle, \\
\langle v_1, v_2, v_3\rangle, &
\langle v_4, v_5, v_6\rangle, & 
\text{and their faces}
\end{array}
\right \}. $$ 
We recommend the reader to draw pictures of 
$\Delta _P$ and $\Delta_Q$ by himself. 
\begin{lem}
We have the following properties{\em{:}}  
\begin{itemize}
\item[(i)] $X_P:=X(\Delta_P)$ is a non-projective 
complete toric variety with $\rho(X_P)=1$, 
\item[(ii)] $X_Q:=X(\Delta _Q)$ is a projective toric 
variety with $\rho(X_Q)=1$, 
\item[(iii)] there exists a toric birational morphism $f_{PQ}: 
X_P\longrightarrow X_Q$, which contracts a $\mathbb P^1$ 
on $X_P$, 
\item[(iv)] $X_P$ and $X_Q$ have only canonical Gorenstein 
singularities, and 
\item[(v)] 
$$
\begin{matrix}
N_1(X_P)&\simeq &N_1(X_Q)\\
\cup& &\cup\\
{NE}(X_P)&\simeq &{NE}(X_Q). 
\end{matrix}
$$ 
In particular, $NE(X_P)=\overline{NE}(X_P)$ is a half line. 
\end{itemize}
\end{lem}

\begin{proof} 
It is easy to check that 
$X_Q$ is projective and $\rho(X_Q)=1$ 
(cf.~\cite[Theorem 3.2, Example 3.5]{eik}). 
Assume that $X_P$ is projective. 
Then there exists a strict upper convex 
support function $h$. 
We note that 
\begin{eqnarray*}
v_1+v_5&=&v_2+v_4,\\
v_2+v_6&=&v_3+v_5, \\
v_3+v_4&=&v_1+v_6. 
\end{eqnarray*}
Thus, we obtain 
\begin{eqnarray*}
h(v_1)+h(v_5)&<&h(v_2)+h(v_4),\\
h(v_2)+h(v_6)&=&h(v_3)+h(v_5), \\
h(v_3)+h(v_4)&=&h(v_1)+h(v_6). 
\end{eqnarray*}
This implies that 
$$
\sum _{i=1}^{6}h(v_i)<\sum _{i=1}^{6}h(v_i). 
$$
It is a contradiction. 
Therefore, $X_P$ is not projective. 
Thus $f_{PQ}$ is not a projective morphism. 
So, $L\cdot C=0$ for 
every $L\in \Pic (X_P)$, where $C\simeq \mathbb P^1$ is 
the exceptional locus of $f_{PQ}$. We note that 
the condition 
(ii) (b) in \cite[p.325 Theorem 1]{kl} does not hold. 
The other statements are easy to check. 
\end{proof}
Let $H$ be an ample Cartier divisor on $X_Q$ 
and $D:=(f_{PQ})^*H$. 
Then $D$ is positive on ${\overline {NE}}(X_P)\setminus \{0\}
=NE(X_P)\setminus \{0\}$. Thus, 
the interior of the nef cone of $X_P$ is non-empty 
(cf.~\cite[p.327 Proposition 2]{kl}). However, 
$D$ is 
not ample on $X_P$. 
Therefore, Kleiman's ampleness criterion does not hold 
for $X_P$. Note that $X_P$ is not projective nor 
quasi-divisorial in the sense of Kleiman 
(see \cite[p.326 Definition 4]{kl}). 
\end{say}

\begin{cor}
We put  $X:=X_P\times \mathbb P^1\times \cdots\times \mathbb P^1$. 
Then we obtain complete non-projective singular toric 
varieties with $\dim X\geq 4$ for which Kleiman's ampleness criterion 
does not hold. 
Since every complete toric surface is 
$\mathbb Q$-factorial and projective, 
Kleiman's ampleness criterion always holds for toric 
surfaces. 
\end{cor}

\begin{rem}\label{1da} 
Let $X$ be a complete 
$\mathbb Q$-factorial 
algebraic variety. Then 
it is not difficult to see that $X$ is projective 
if $\rho(X)=1$. 
\end{rem}

\begin{rem}
In \cite[Chapter VI. Appendix 2.19.3 Exercise]{ko}, 
Koll\'ar pointed out that Kleiman's ampleness criterion does not hold 
for smooth proper algebraic spaces. 
\end{rem}

\begin{rem}
In \cite[Theorem 4.1]{part1}, we claim that 
$NE(X/Y)$ is strongly convex if $f:X\longrightarrow Y$ is 
projective. This is obvious. 
However, in the proof of Theorem 4.1 in \cite{part1}, we say that 
it follows from Kleiman's criterion. Sorry, it is misleading. 
\end{rem}

We note the following ampleness criterion, which works for 
complete toric varieties with {\em{arbitrary}} singularities. 

\begin{prop}
Let $X$ be a complete toric variety and $L$ a 
line bundle on $X$. 
Assume that $L\cdot C>0$ for every torus invariant 
integral curve $C$ on $X$. 
Then $L$ is ample. 
In particular, $X$ is projective. 
\end{prop}
\begin{proof} 
Since $NE(X)$ is spanned by the torus invariant 
curves on $X$, it is obvious that $L$ is nef. 
This implies that $L$ is generated by its global 
sections. Note that we can replace 
$X$ (resp.~$L$) with its toric resolution $Y$ (resp.~the 
pull-back of $L$ on $Y$) in order to check the freeness of $L$. 
Thus, the proof of 
the freeness is easy. We consider the equivariant 
morphism $f:=\Phi_{|L|}:X\longrightarrow Y$ associated to 
the linear system $|L|$. Then 
we obtain that $L=f^*H$ for a very ample line bundle 
$H$ on $Y$. Since $L\cdot C>0$ for every 
torus invariant integral 
curve $C$ on $X$, we have that $f$ is finite. 
Thus, $L$ is ample.  
\end{proof}

\section{Singular varieties with $NE(X)=N_1(X)$}\label{sec3} 

In this section, we construct complete singular toric 
varieties with $NE(X)=N_1(X)\simeq \mathbb R^k$ for any $k\geq 0$. 

\begin{rem}
The condition $NE(X)=N_1(X)$ is 
equivalent to the following one:~a line 
bundle $L$ is nef if and only if $L$ is numerically equivalent to 
zero. 
\end{rem}

\begin{say}[Construction]\label{31} 
We fix $N=\mathbb Z^3$ and $M:=\Hom_{\mathbb Z}(N,\mathbb Z)\simeq 
\mathbb Z^3$. 
We put 
\begin{align*}
v_1  &= (1,0,1),  &v_2&=(0,1,1),  &v_3&=(-1,-2,1),\\
v_4  &= (1,0,-1),  &v_5&=(0,1,-1),  &v_6&=(-1,-1,-1), \\ 
v_7  &= (0,0,-1).  & & &
\end{align*}
First, we consider the following fan. 
$$
\Delta_A=
\left\{ 
\begin{array}{ccc}
\langle v_1, v_2, v_4, v_5\rangle, &
\langle v_2, v_3, v_5, v_6\rangle, &
\langle v_1, v_3, v_4, v_6\rangle, \\
\langle v_1, v_2, v_3\rangle, &
\langle v_4, v_5, v_6\rangle, & 
\text{and their faces}
\end{array}
\right \}. $$
We recommend the reader to draw the picture of 
$\Delta _A$ by himself. 
We put $X_A:=X(\Delta_A)$ and $D_i:=V(v_i)$ for 
every $i$. This example $X_A$ is essentially the same as 
\cite[Example 3.5]{eik}.  
We consider the principal divisor 
$$
D=D_1+D_2+D_3-D_4-D_5-D_6
$$ 
that is associated to $m=(0, 0, -1)\in M$. 
We put 
\begin{align*}
\sigma_1&=\langle v_1, v_2, v_4, v_5\rangle, &
\sigma_2&=\langle v_2, v_3, v_5, v_6\rangle, &
\sigma_3&=\langle v_1, v_3, v_4, v_6\rangle. \\ 
\end{align*}
Then, all points in $\sum _i\AD(\sigma_i)$ are linear 
combinations of the lines of the matrix 
$$
\begin{pmatrix}
1&-1&0&1&-1&0\\
0&-1&2&0&-3&2\\
-2&0&1&-1&0&2\\
\end{pmatrix}
$$
which has rank $3$. 
Note that we use $D$ to define $\AD(\sigma_i)$ 
and 
that $\sigma_i$ is not simplicial for every $i$ and 
all the other $3$-dimensional cones in $\Delta_A$ are simplicial. 
For the definition of $\AD(\sigma_i)$, see \cite[Theorem 3.2]{eik}. 
Therefore, $\Pic X_A\simeq \mathbb Z^{6-3-3}=\{0\}$ by 
\cite[Theorem 3.2]{eik}. 

Next, we consider the following fan. 
$$
\Delta_B=
\left\{ 
\begin{array}{ccc}
\langle v_1, v_2, v_4, v_5\rangle, &
\langle v_2, v_3, v_5, v_6\rangle, &
\langle v_1, v_3, v_4, v_6\rangle, \\
\langle v_1, v_2, v_3\rangle, &
\langle v_4, v_5, v_6\rangle, &
\langle v_4, v_5, v_7\rangle, \\
\langle v_4, v_6, v_7\rangle, &
\langle v_5, v_6, v_7\rangle, & 
\text{and their faces}
\end{array}
\right \}. $$ 
We recommend the reader to draw the picture of 
$\Delta _B$ by himself. 
We put $X_B:=X(\Delta_B)$. 
Then $X_B\longrightarrow X_A$ is the blow up along 
$\langle v_7\rangle$. 
We consider the principal divisor 
$$
D'=D_1+D_2+D_3-D_4-D_5-D_6-D_7
$$ 
that is associated to $m=(0, 0, -1)\in M$. 
Then, all points in $\sum_{i=1}^{3}\AD(\sigma_i)$ are linear 
combinations of the lines of the 
matrix 
$$
\begin{pmatrix}
1&-1&0&1&-1&0&0\\
0&-1&2&0&-3&2&0\\
-2&0&1&-1&0&2&0\\
\end{pmatrix}
$$
which has rank $3$. 
We note that we use $D'$ to define $\AD(\sigma_i)$ and 
that $\sigma_i$ is not simplicial for every $i$ and 
all the other $3$-dimensional cones in $\Delta_B$ are simplicial. 
Thus, we have $\Pic X_B\simeq \mathbb Z^{7-3-3}= \mathbb Z$ 
by \cite[Theorem 3.2]{eik}. 
\begin{lem}
$D_7$ is a Cartier divisor. 
\end{lem}
\begin{proof}
It is because each $3$-dimensional 
cone containing $v_7$ is non-singular. 
Thus, $D_7$ is Cartier. 
\end{proof}
We put $C_1:=V(\langle v_4, v_5\rangle)\simeq \mathbb P^1$ and 
$C_2:=V(\langle v_4, v_7\rangle)\simeq \mathbb P^1$. 
The following lemma is a key property of this example. 
\begin{lem}
$C_1\cdot D_7>0$ and $C_2\cdot D_7<0$. 
More precisely, $C_1\cdot D_7=1$ and $C_2\cdot D_7=-3$. 
Therefore, $D_7$ is a generator of $\Pic X_B\simeq \mathbb Z$.  
\end{lem}
\begin{proof} 
It is obvious that $C_1\cdot  D_7=1>0$. 
Since $v_4+v_5+v_6-3v_7=0$, 
we have $C_2\cdot D_7=-3C_2\cdot D_5=-3<0$. 
\end{proof} 
Therefore, $NE(X_B)=N_1(X_B)\simeq \mathbb R$. 
In particular, $X_B$ is not projective. 
\end{say}
\begin{cor}
Let $X:=X_A$. 
Then $NE(X)=N_1(X)=\{0\}$. 
Let $X:=X_B\times X_B\times \cdots\times X_B$ be the $k$ times 
product of $X_B$. 
Then $NE(X)=N_1(X)\simeq \mathbb R^k$. 
\end{cor}
\begin{proof}
The first statement is obvious by the above construction. 
It is not difficult to see that 
$\Pic X\simeq \bigotimes_{i=1}^{k}p_{i}^{*}\Pic X_B$, 
where $p_i:X\longrightarrow X_B$ is the $i$-th projection. Thus, 
we can check that any nef line bundle on $X$ is 
numerically trivial. 
Therefore, $NE(X)=N_1(X)\simeq \mathbb R^k$. 
\end{proof}

In the above corollary, $\dim X=3k$ if $NE(X)=N_1(X)
\simeq \mathbb R^k$ for $k\geq 1$. We can 
construct $3$-folds with $NE(X)=N_1(X)\simeq 
\mathbb R^k$ for every $k\geq 0$. 
We just note this fact in the next remark. 
The details are left to the reader. 

\begin{rem} 
We put $X_0:=X_A$, where $X_A$ is in \ref{31}. 
We define primitive vectors $\{v_k\}$ inductively 
for $k\geq 7$ as follows:~
$\langle v_k\rangle =\langle v_4+v_5+v_{k-1}\rangle$. 
Let $X_k\longrightarrow X_{k-1}$ be the blow up 
along $\langle v_{k+6}\rangle$ for $k\geq 1$. 
Note that $X_1=X_B$, where $X_B$ is in \ref{31}. 
It is not difficult to see that 
$NE(X_k)=N_1(X_k)\simeq \mathbb R^k$ and 
$\dim X_k=3$. 
We note that the numerical equivalence classes 
of $D_i=V(v_i)$'s for $7\leq i\leq k+6$ form 
a basis of $N^1(X_k)$.  
\end{rem}

\section{Miscellaneous comments}

\begin{say}
Let $X$ be a complete normal 
variety. Assume that $X$ is $\mathbb Q$-factorial. 
We note 
that $\rho(X)\geq 1$ in this case. 
Our question is 
as follows. 

\begin{que}\label{444} 
Are there any complete normal $\mathbb Q$-factorial 
varieties $X$ with 
$NE(X)=N_1(X)$? 
\end{que} 
The examples constructed in Section \ref{sec3} are obviously 
non-$\mathbb Q$-factorial. 
We note that $\rho(X)\geq 2$ when $NE(X)=N_1(X)$ by 
Remark \ref{1da}. 

\begin{rem} 
Let $X$ be a complete normal (not necessarily 
$\mathbb Q$-factorial) variety. 
If there exists a proper surjective 
morphism $f:X\longrightarrow Y$ such 
that 
\begin{itemize}
\item[(i)] $f$ has connected fibers, 
\item[(ii)] $Y$ is projective and $\dim Y\geq 1$, and 
\item[(iii)] $f$ is not an isomorphism, 
\end{itemize}
then it is obvious that $NE(X)\subset \overline{NE}(X)
\subsetneq N_1(X)$. 
\end{rem}
\end{say}

\begin{say}[Toric variety] 
We have the following conjecture on Question \ref{444}.  

\begin{conj}
Let $X$ be a non-singular complete toric 
variety. 
Then $\overline{NE}(X)=NE(X)\subsetneq N_1(X)$. 
\end{conj}

\begin{rem}
If $X$ is a complete (not necessarily $\mathbb Q$-factorial) 
toric variety with $NE(X)\subsetneq 
N_1(X)$, then there exists a non-trivial nef line bundle 
$\mathcal L$ on $X$. Thus, we obtain the toric morphism 
$f:X\longrightarrow Y:=\Proj 
\bigoplus _{k\geq 0}H^0(X, \mathcal 
L^{\otimes k})$ such that $\mathcal L=f^*\mathcal H$, 
where $\mathcal H$ is an ample line bundle on $Y$. 
Note that $\dim Y\geq 1$. 
\end{rem}
\end{say}

\begin{say}[Algebraic space]  
I learned the following example from S.~Mori, who call 
it Hironaka's example. I have never seen this in 
the literature. 

\begin{ex}
Let $Q\simeq \mathbb P^1\times \mathbb P^1$ 
be a non-singular quadric surface in $\mathbb P^3_{\mathbb 
C}$. 
We take a non-singular $(3,d)$-curve $C$ in $Q$, 
where $d\in \mathbb Z_{>0}$. 
That is, $\mathcal O_Q(C)\simeq p_1^*\mathcal O_{\mathbb P^1}(3)
\otimes p_2^*\mathcal O_{\mathbb P^1}(d)$, 
where $p_1$ (resp.~$p_2$) is the first (resp.~second) projection from 
$Q$ to $\mathbb P^1$. 
Let $f_1$ (resp.~$f_2$) be a fiber of $p_2:Q\longrightarrow \mathbb P^1$ 
(resp.~$p_1$). 
We take the blow up $\pi:X\longrightarrow \mathbb P^3$ along $C$. 
Let $Q'$ be the strict transform of $Q$ and $E$ the exceptional divisor 
of $\pi$. 
Then $Q'=f^*Q-E$. 
Let $f_i'$ be the strict transform of $f_i$ for $i=1,2$. 
We have 
$$
Q'\cdot f_1'=f^*Q\cdot f_1'-E\cdot f_1'=Q\cdot f_1-3=2-3=-1. 
$$ 
Thus we can blow down $X$ to $Y$ along the ruling $p_2:Q'\simeq 
Q\longrightarrow \mathbb P^1$ 
(cf.~\cite[Main Theorem]{nakano}, \cite{fn}). 
Note that $Y$ is a compact Moishezon manifold. 
The Kleiman-Mori cone $\overline {NE}(X)$ is spanned by $2$ rays $R$ and 
$Q$. We note that $X$ is non-singular projective and 
$\rho(X)=2$. 
Let $l$ be a fiber of $\pi:X\longrightarrow \mathbb P^3$. 
Then, one ray $R$ is spanned by the numerical equivalence 
class of $l$. 
We put $\mathcal L:=\pi^*\mathcal O_{\mathbb P^3}(1)$. 
Then $\mathcal L$ is non-negative on $\overline{NE}(X)$ 
and $R=(\mathcal L=0)\cap \overline {NE}(X)$. 
We have the following intersection numbers. 
$$
\mathcal L\cdot l=0, \ \ 
\mathcal L\cdot f_1'=\mathcal L\cdot f_2'=2, 
$$
$$
E\cdot l=-1, \ \ 
E\cdot f_1'=3, \ \ 
E\cdot f_2'=d. 
$$
From now on, we assume $d\geq 4$. 
We can write $f_1'=af_2'+bl$ in $N_1(X)$ for $a, b \in \mathbb R$. 
Thus we can easily check that $a=1$, $b=d-3>0$. 
Therefore, the numerical class of $f_1'$ 
is in the interior of the cone 
spanned by the numerical classes of $f_2'$ and $l$. 
Thus, we have 
$NE(Y)=\overline {NE}(Y)=N_1(Y)\simeq\mathbb R$. 
Therefore, $Y$ is a non-singular complete algebraic space with 
$\rho(X)=1$. Note that $Y$ is not a scheme. 
\end{ex}
\end{say} 
\ifx\undefined\bysame
\newcommand{\bysame|{leavemode\hbox to3em{\hrulefill}\,}
\fi

\end{document}